\documentclass[12pt]{amsart}
\usepackage{amsmath}
\usepackage{amssymb}
\def\Ind#1{1_{\{#1\}}}
\def\Cal{\mathcal}

\def\e{\varepsilon}

\def\oo{\infty}

\def\half{{\tfrac12}}

\def\upref#1{\upn{\ref{#1}}}

\theoremstyle{definition}

\theoremstyle{plain}
\newtheorem{thm}{Theorem}

\newtheorem{lem}{Lemma}

\theoremstyle{remark}
\newtheorem{rmrk}{Remark}

\title[A comparison inequality]{A comparison inequality for sums of independent
random variables}
\author{Stephen J. Montgomery-Smith}
\author{Alexander R. Pruss}
\date{November 20, 1998}
\address{Department of Mathematics\\
University of Missouri\\
Columbia, MO 65211\\
U.S.A.}
\address{Department of Philosophy\\
University of Pittsburgh\\
Pittsburgh, PA 15260\\
U.S.A.}
\email{stephen@math.missouri.edu}
\email{pruss+@pitt.edu}
\urladdr{http://www.missouri.edu/\~{}stephen}
\urladdr{http://www.pitt.edu/\~{}pruss}
\begin{document}
\textwidth=30cc
\baselineskip=16pt
\begin{abstract}
        We give a comparison inequality that allows one to estimate the
        tail probabilities of sums of independent
        Banach space valued
        random variables in
        terms of those of independent {\em identically distributed}
        random variables.
        More precisely, let $X_1,\dots,X_n$ be independent
        Banach-valued
        random variables.  Let $I$ be a random variable independent of
        $X_1,\dots,X_n$ and uniformly distributed over $\{ 1,\dots,n
        \}$.  Put $\tilde X_1=X_I$, and let $\tilde X_2,\dots,\tilde
        X_n$ be independent identically distributed copies of $\tilde
        X_1$.
         Then,
        $P(\|X_1+\dots+X_n\| \ge \lambda)\le cP(\|\tilde
        X_1+\dots+\tilde X_n\|\ge \lambda/c)$ for all
        $\lambda\ge 0$, where $c$ is an absolute constant.
\end{abstract}

\maketitle
\bibliographystyle{amsplain}
The independent Banach-valued
random
variables $X_1,\dots,X_n$ are said
to {\em regularly cover} (the distribution of) a random variable $Y$
provided that
$$
        E[g(Y)]=\frac1n \sum_{k=1}^n E[g(X_k)],
$$
for all Borel functions $g$ for which either side is
defined~\cite{Pruss:Riemann}.  An easy way of constructing $Y$, given
the independent Banach-valued random variables $X_1,\dots,X_n$, is to let
$I$ be a random variable independent of $X_1,\dots,X_n$, with values in
$\{1,2,\dots,n\}$ and with each value having equal probability $1/n$,
and then put $Y=X_I$.  It is easy to see that then $X_1,\dots,X_n$
regularly cover $Y$.  This construction will be useful for our proofs.

If the variables are real valued, then the regular covering condition is
easily seen to be equivalent to the condition that the distribution
function $F$ of $Y$ is the arithmetic mean of the respective
distribution functions $F_1,\dots,F_n$ of $X_1,\dots,X_n$.

A variable $X'$ is said to be a {\em copy}\/ of $X$ if it has the same
distribution as $X$. The main purpose of this paper is then to prove the
following result.

\begin{thm}\label{th:main}
        There exists an absolute constant $c\in (0,\oo)$ such that if
        $X_1,\dots,X_n$ are independent Banach-valued random variables which
regularly
        cover a random variable $\tilde X_1$, then:
\begin{equation}\label{eq:the-ineq}
        P(\|X_1+\dots+X_n\|\ge\lambda) \le cP(\|\tilde X_1+\dots+\tilde
        X_n\|\ge\lambda/c),
\end{equation}
        for all $\lambda\ge 0$, where $\tilde X_2,\dots,\tilde X_n$
        are independent copies of $\tilde X_1$.
\end{thm}

\begin{rmrk}
        In the case where the random variables are symmetric, this was
        shown in \cite{Pruss:noniid} (strictly speaking, it was only
        shown in the real-valued case, but the proof also works for the
        Banach-valued case).
\end{rmrk}

\begin{rmrk}
        The inequality converse to \eqref{eq:the-ineq} is false, even in
        the special cases of symmetric real random variables.
        For, suppose that $c$ is an absolute constant such that
\begin{equation}\label{eq:conv-ineq}
        P(|\tilde X_1+\dots+\tilde X_n|\ge\lambda) \le cP(|X_1+\dots+
        X_n|\ge\lambda/c),
\end{equation}
        for all $\lambda\ge 0$, whenever the conditions of
        Theorem~\ref{th:main} hold with symmetric variables. Fix any
        $n>\max(1,c)$. Put
        $X_2\equiv\dots\equiv X_n\equiv 0$.  Let $X_1$ be such that
        $P(X_1=1)=P(X_1=-1)=\half$.  Put $\lambda=n$.  Then the right
        hand side of \eqref{eq:conv-ineq} is zero, since
        $|X_1+\dots+X_n|\equiv 1$.  But the left hand side of
        \eqref{eq:conv-ineq} is non-zero, since it is easy to see that
        $P(\tilde X_i=1)=2^{-n-1}$ for each $i$ (as the $\tilde X_i$ are
        identically distributed, and as $\tilde X_1$ can be taken to be
        $X_I$ where $I$ is independent of everything else and uniformly
        distributed on $\{1,\dots,n\}$), so that $P(|\tilde X_1+\dots+\tilde
        X_n|\ge n)\ge (2^{-n-1})^n > 0$.
\end{rmrk}

\begin{rmrk}
The main consequence of Theorem~\upref{th:main} is that any upper
bound on tail probabilities of sums of independent identically
distributed random variables automatically gives a bound on tail
probabilities of sums of non-identically distributed independent
random variables.
\end{rmrk}

\begin{rmrk}
For a very simple application, we give another proof of one side of a
result from \cite{Pruss:Riemann} on randomly sampled Riemann sums.  Let $f\in L^2[0,1]$.
For $1\le k\le n$, let $x_{nk}$ be
uniformly distributed over $[(k-1)/n,k/n]$, and assume
$x_{n1},\dots,x_{nn}$ are independent for each fixed $n$.  Define the randomly sampled
Riemann sum $R_n f=n^{-1}\sum_{k=1}^n f(x_{nk})$.  Then the result says
that $R_n f$ converges almost surely to the Lebesgue integral
$A=\int_0^1 f$.  (For a converse in the case where {\em all}\/ the
$x_{nk}$ are independent, not just for fixed $n$, see
\cite{Pruss:Riemann}.)  For, by Borel-Cantelli it suffices to show that
\begin{equation}\label{eq:rsum}
        \sum_{n=1}^\oo P(|R_n f-A|\ge \e)<\oo,
\end{equation}
for all $\e>0$. Let $X_1,X_2,\dots$ be independent identically distributed
random variables with the same distribution as $f$. Note that
$f(x_{n1}),\dots,f(x_{nn})$ regularly cover $X_1$, and
$f(x_{n1})-A,\dots,f(x_{nn})-A$ regularly cover $X_1-A$.  Since $f\in
L^2$, we have $X_1$ having a finite second moment, and moreover
$E[X_1]=A$, so that by the
Hsu-Robbins law of large numbers~\cite{HsuRobbins} (see also
\cite{Erdos,Erdos:Two}), we have
$$
        \sum_{n=1}^\oo P(|(X_1-A)+\dots+(X_n-A)|/n \ge \e)<\oo,
$$
for all $\e>0$. By Theorem~\ref{th:main} and the fact that
$f(x_{n1})-A,\dots,f(x_{n1})-A$ regularly cover $X_1-A$, we obtain
\eqref{eq:rsum}.
\end{rmrk}

To prove Theorem~\upref{th:main}, we need some definitions and lemmata.
If $X$ is a random variable, then let $X^s=X-X'$ be the {\em
symmetrization}\/ of $X$, where $X'$ is an independent copy of $X$.  We
shall always choose symmetrizations so that we have
$(X_1+\dots+X_k)^s=X_1^s+\dots+X_k^s$ whenever we need this identity.

Write $\|X\|_p = (E[\|X\|^p])^{1/p}$, where $\| \cdot \|$ is the norm on
the Banach space in which our random variables take values.

\begin{lem}\label{lem:disymm2}
        Let $X$ be a Banach-valued 
        random variable with $\| X \|_2<\oo$.
        Then, $\| X \|_2 \le \|X^s\|_2 + \|E [X]\| \le 3 \| X \|_2$
\end{lem}

\begin{proof}
Let $X'$ be an independent copy of $X$ so that $X^s=X-X'$.
Let $\Cal A$ be the sigma-algebra generated by $X$.
Then $E[X^s\mid\Cal A] = X-E[X']=X-E[X]$, and so
\[
\|X\|_2 = \| E[X^s+E[X]\mid \Cal A] \|_2
  \le \| X^s+E[X] \|_2
  \le \| X^s \|_2 + \| E[X] \|,
\]
where the first inequality used the fact that conditional expectation is
a contraction on the Banach-valued $L^p$ spaces, $p\ge 1$ (see, e.g.,
\cite[Theorem~V.1.4]{DiestelUhl}).
The rest of the Lemma follows from the triangle inequality.
\end{proof}

\begin{lem}\label{lem:comp-moment}
Let $X_1,\dots,X_n$ be independent random variables, and let
$\tilde X_1,\dots,\tilde X_n$ be independent identically distributed
random variables such that $X_1,\dots,X_n$ regularly cover $\tilde X_1$.
Put $S_n=X_1+\dots+X_n$ and $\tilde S_n=\tilde X_1+\dots+\tilde X_n$.
Then:
\[
         \| S_n \|_2 \le 12  \| \tilde S_n \|_2.
\]

\end{lem}

\begin{proof}  Let $I_1,\dots,I_n$ be independent random variables
uniformly distributed on the set $\{ 1,\dots,n \}$.  Let $\{ X_{i,j}
\}_{1\le i,j\le n}$ and $\{ X_{i,j}' \}_{1\le i,j\le n}$ be independent
arrays of independent random variables, with the arrays independent of
the $I_i$, and such that $X_{i,j}$ and $X_{i,j}'$ both have the same
distribution as $X_j$ for all $i$ and $j$.  Without loss of generality
we can put $\tilde X_i=X_{i,I_i}$.  Set $\tilde X_i'=X'_{i,I_i}$.
Let $\tilde S_n'=\tilde X_1'+\dots+\tilde X_n'$.  Let
$(X_1',\dots,X_n')$ be an independent copy of $(X_1,\dots,X_n)$, and put
$S_n'=X_1'+\dots+X_n'$.  Observe that $X_1-X_1',\dots,X_n-X_n'$
regularly cover $\tilde X_i-\tilde X_i'$ for all $i$, and that moreover
the $X_i-X_i'$ are symmetric.  Thus, by
\cite[Proposition~1]{Pruss:noniid} (which though stated for real valued
random variables, holds for the Banach-valued
case as well, and with the same proof) we have:
\begin{equation}\label{eq:moment-1}
        \| S_n-S_n' \|_2 \le 4 \| \tilde S_n-\tilde S_n' \|_2.
\end{equation}
Also, it is clear that $E[S_n] = E[\tilde S_n]$.
Combining this with Lemma~\upref{lem:disymm2}, we see that:
\[
     \| S_n \|_2 \le \| S_n-S_n' \|_2 + \| E[S_n] \|
     \le 4 \| \tilde S_n -\tilde S_n'\|_2 + 4 \|
     E[\tilde S_n] \|
     \le 12 \| \tilde S_n \|_2 ,
\]
as desired.
\end{proof}

The following Lemma is in effect a special case of a result of
Hitczenko~\cite{Hitczenko}.

\begin{lem}\label{lem:Rosenthal}
        Let $X_1,\dots,X_n$ be independent identically distributed 
	Banach-valued
        random variables with $\| X_i \|<L$ almost surely for all
        $i$.  Let $S_k=X_1+\dots+X_k$.  Then:
\[
        (E[\| S_n \|])^2 \ge c (E[ \| S_n \|^2 ]-c^{-1} L^2),
\]
        where $c\in (0,\oo)$ is an absolute constants.
\end{lem}

\begin{proof}[Proof of Lemma~\upref{lem:Rosenthal}]
By the work of Hitczenko~\cite{Hitczenko},
if $S^* = \max_{k}\|S_k\|$ and
$X^* = \max_k\|X_k\|$, then for $q \ge p$:
$$ \|S^*\|_q \le c_0 \frac qp ( \|S^*\|_p + \|X^*\|_q ),$$
for a finite absolute constant $c_0$.
By \cite[Corollary~4]{MS} we have $\| S^* \|_p \le c_1 \| S_n \|_p$ for an
absolute constant $c_1$, as the $X_i$ are identically distributed.  The
desired inequality easily follows from this
with $c=8c_0^2$ if we let $q=2$ and $p=1$.
\end{proof}

\begin{proof}[Proof of Theorem~\upref{th:main}]
Let $I_1,\dots,I_n$, $\{ X_{i,j} \}_{1\le i,j\le n}$, $\{
X_{i,j}' \}_{1\le i,j\le n}$, $S_n'$ and $\tilde S_n'$ be as in the
proof of Lemma~\upref{lem:comp-moment}. Applying
\cite[Proposition~1]{Pruss:noniid} (which works for
Banach-valued variables as already stated), we see that
\begin{equation}\label{eq:first-ineq}
   P(\| S_n-S_n' \|\ge \lambda) \le 8P(\| \tilde S_n-\tilde S_n' \| \ge
   \lambda/2) \le 16P(\| \tilde S_n \| \ge
   \lambda/4),
\end{equation}
for all $\lambda$,
where the second inequality followed from the
inequality that $P(\|X^s\|\ge t)\le P(\|X\|\ge t/2)+P(\|X'\|\ge
t/2)=2P(\|X\|\ge
t/2)$, where $X'$ is an independent copy of $X$ such that $X^s=X-X'$.
Note that $S_n^s=S_n-S_n'$.

Let $M$ be a median of $\| S_n \|$.  It is easy to see that
\begin{equation}\label{eq:second-ineq}
        P(\|S_n\|-M\ge\lambda)\le 2P(|S_n^s|\ge\lambda),
\end{equation}
for all $\lambda$.  (For, if $\| S_n \|-M\ge\lambda$, there is at least
probability $1/2$ that $\| S_n' \|\le M$ in which case $\| S_n -S_n'\|
\ge \| S_n \|-\|S_n'\| \ge \|S_n \|-M\ge\lambda$.)

We now claim that in general in our present setting:
\begin{equation}\label{eq:claim}
        P(\|\tilde S_n\|\ge \e M)\ge \delta,
\end{equation}
for absolute constants $\e,\delta\in (0,1)$ to be determined later.
(they will be determined in accordance with
\eqref{eq:1a-e-ass},
\eqref{eq:second-e-ass},
\eqref{eq:first-e-ass},
\eqref{eq:third-e-ass} and
\eqref{eq:fourth-e-ass}, below).  To prove \eqref{eq:claim}, suppose
that on the contrary we have:
\begin{equation}\label{eq:n-claim}
        P(\|\tilde S_n\|\ge \e M)\le \delta.
\end{equation}
Since the $\tilde X_i$ are independent and identically distributed,
by a maximal inequality for sums of independent and identically
distributed random variables~\cite[Corollary~4]{MS} together with
\eqref{eq:n-claim}, we have:
\begin{equation}\label{eq:star1}
        P\bigl(\max_{1\le k\le n} \|\tilde S_k\|\ge c_1\e M\bigr)
         \le c_1P(\|\tilde S_n\|\ge \e M)\le c_1\delta,
\end{equation}
where $c_1\in [1,\oo)$ is an absolute constant.  By the elementary
inequality
$$
        P(\max_{1\le k\le n} \|U_k\|\ge 2t)\le
        P\biggl(\max_{1\le k\le n} \biggl\| \sum_{i=1}^k U_i \biggr\| \ge
        t\biggr),
$$
valid for all $t$ if the $U_i$ are independent (since
if $\|U_k\|\ge 2t$ then $\| \sum_{i=1}^k U_i \|\ge t$ or $\|\sum_{i=1}^{k-1}
U_i\|\ge t$), it follows from \eqref{eq:star1} that
\begin{equation}\label{eq:smallness}
        P\bigl(\max_{1\le k\le n} \|\tilde X_k\|\ge 2c_1\e M\bigr)
         \le c_1\delta.
\end{equation}
Let $L=2c_1\e M$.
Set $Y_k=X_k\cdot \Ind{ \|X_k\|<L }$.  Put $\tilde Y_k=\tilde
X_k\cdot\Ind{\|\tilde X_k\|< L }$.  Note that $Y_1,\dots,Y_n$ regularly
cover $\tilde Y_k$ for each $k$.  Let $T_n=Y_1+\dots+Y_n$ and put
$\tilde T_n=\tilde Y_1+\dots+\tilde Y_n$.  By \eqref{eq:smallness}, we
have:
\begin{equation}\label{eq:tilde-close}
        P\biggl(\bigcup_{k=1}^n \{ \tilde X_k\ne\tilde Y_k \}\biggr)
         \le c_1\delta.
\end{equation}
Let $p=P(\|\tilde X_k\|\ge L)$.  Note that this does not depend on
$k$ since the $\tilde X_k$ are identically distributed.  Note also that
the left hand side of \eqref{eq:tilde-close} is equal to $1-(1-p)^n$.
Henceforth we will assume that
\begin{equation}\label{eq:1a-e-ass}
        \delta<1/(2c_1).
\end{equation}
Now, if $x\in[0,1]$ is such that $1-(1-x)^n\le 1/2$, then $nx\le 2(1-(1-x)^n)$.
Then, using this observation, together with
\eqref{eq:tilde-close}, \eqref{eq:1a-e-ass}
and the condition that $X_1,\dots,X_n$ regularly cover $\tilde X_1$:
\begin{equation}\label{eq:nontilde-close}
\begin{split}
        P\biggl(\bigcup_{k=1}^n \{ X_k\ne Y_k \}\biggr)
         &\le \sum_{k=1}^n P(X_k\ne Y_k) \\
         &=\sum_{k=1}^n P(\|X_k\|\ge L) \\
         &=n P(\|\tilde X_1\|\ge L)
         =np\\
         &\le 2(1-(1-p))^n \\
         &= 2P\biggl(\bigcup_{k=1}^n \{ \tilde X_k\ne\tilde Y_k \}\biggr)
         \le 2c_1\delta.
\end{split}
\end{equation}

Now, by \eqref{eq:first-ineq}, \eqref{eq:second-ineq} and
\eqref{eq:n-claim}, it follows that
$$
        P(\|S_n\|-M\ge 4\e M) \le 32\delta.
$$
Using \eqref{eq:nontilde-close}, it then follows that:
\begin{equation}\label{eq:n-claim-2-adj}
        P(\|T_n\|-M\ge 4\e M) \le (32+2c_1)\delta.
\end{equation}
Moreover, by \eqref{eq:n-claim} and \eqref{eq:tilde-close}:
\begin{equation}\label{eq:n-claim-s}
        P(\|\tilde T_n\|\ge \e M)\le (1+c_1)\delta.
\end{equation}

Observe that $|\tilde Y_i|<L$ almost surely.
Lemma~\ref{lem:Rosenthal} then shows that:
\begin{equation}\label{eq:Rosenthal}
        (E[ \|\tilde T_n\| ])^2 \ge c_2(E[\|\tilde T_n\|^2]-c_2^{-1} L^2 ),
\end{equation}
where $c_2\in (0,\oo)$ is an absolute constant.

Now, by \eqref{eq:n-claim-2-adj} we have:
\begin{equation}\label{eq:lower-est-22}
   E[\|T_n\|^2]\ge [1-(32+2c_1)\delta]M^2.
\end{equation}
Henceforth, we will assume that $\delta$ is sufficiently small that
\begin{equation}\label{eq:second-e-ass}
        1-(32+2c_1)\delta\ge \half.
\end{equation}
Using Lemma~\upref{lem:comp-moment} we see that
$E[\|T_n\|^2] \le 144 E[\| \tilde T_n \|^2]$.  Combining this with
\eqref{eq:lower-est-22} and \eqref{eq:second-e-ass}, we see that
\begin{equation}\label{eq:2nd-big}
        E[\|\tilde T_n\|^2]\ge M^2/288.
\end{equation}
Assume that $\e>0$ is sufficiently small that $c_2^{-1} L^2\le
M^2/(2\cdot 288)$.
Since $L=2c_1\e M$, this assumption is equivalent to:
\begin{equation}\label{eq:first-e-ass}
        \e\le (48c_1)^{-1}c_2^{1/2}.
\end{equation}
Thus by \eqref{eq:2nd-big}:
\begin{equation}\label{eq:L-small}
  c_2^{-1} L^2 \le E[\|\tilde T_n\|^2]/2.
\end{equation}
Then, by \eqref{eq:Rosenthal},
\begin{equation}\label{eq:mom-ineq}
  (E[\|\tilde T_n\|])^2 \ge c_2 (E[\| \tilde T_n\|^2] - c_2^{-1} L^2)
  \ge \half c_2 E[\| \tilde T_n\|^2].
\end{equation}
The elementary inequality $P(|\Xi|\ge\lambda E[|\Xi|])\ge (1-\lambda)^2
(E[|\Xi|])^2/E[|\Xi|^2]$ (see, e.g., \cite[Exercise 3.3.11]{Chung}) then
implies that
\begin{equation}\label{eq:star2}
  P(\|\tilde T_n\|\ge \half E[\|\tilde T_n\|])
   \ge (1-\half)^2\cdot \half c_2.
\end{equation}
Now, by \eqref{eq:2nd-big} and
\eqref{eq:mom-ineq} we have $E[\|\tilde T_n\|]\ge (\half
c_2/288)^{1/2} M=c_2^{1/2}M/24$, so
that \eqref{eq:star2} gives:
\begin{equation}\label{eq:concl}
  P(\|\tilde T_n\|\ge \half c_2^{1/2} M/24)
    \ge c_2/8.
\end{equation}
If we choose $\e$ and $\delta$ such that
\begin{equation}\label{eq:third-e-ass}
        0<\e\le c_2^{1/2}/48
\end{equation}
and
\begin{equation}\label{eq:fourth-e-ass}
0<(1+c_1)\delta<c_2/8
\end{equation}
and satisfying the
other conditions required in the above argument (namely
\eqref{eq:1a-e-ass}, \eqref{eq:second-e-ass} and \eqref{eq:first-e-ass}),
we will obtain from \eqref{eq:concl} a
contradiction to \eqref{eq:n-claim-s}.
Hence, if we take $\e$ and $\delta$ to be absolute constants in $(0,1)$
satisfying these assumptions, we obtain \eqref{eq:claim}.

Now, combining \eqref{eq:first-ineq} and \eqref{eq:second-ineq}, we see
that:
\begin{equation}\label{eq:third-ineq}
   P(\|S_n\|-M\ge\lambda)\le 32 P(\|\tilde S_n\|\ge\lambda/4),
\end{equation}
for all $\lambda$.
There are now two cases to be considered.  Suppose first
that $\lambda\le 2M$.  Then using \eqref{eq:claim}:
\begin{equation}\label{eq:case-1}
  P(\|S_n\|\ge\lambda)
   \le 1 \le \delta^{-1} P(\|\tilde S_n\|\ge \e M)
   \le \delta^{-1}P(\|\tilde S_n\|\ge \e\lambda/2).
\end{equation}
On the other hand, suppose that $\lambda> 2M$.  In that case
if $\|S_n\|\ge\lambda$ then $\|S_n\| -M> \lambda-\lambda/2=\lambda/2$, so
that
\begin{equation}\label{eq:case-2}
  P(\|S_n\|\ge\lambda)\le P(\|S_n\|-M\ge\lambda/2) \le 32 P(\|\tilde S_n\|\ge
   \lambda/4),
\end{equation}
by \eqref{eq:third-ineq}.  Inequality \eqref{eq:the-ineq} follows from
\eqref{eq:case-1} for $\lambda\le 2M$ and from \eqref{eq:case-2} for
$\lambda > 2M$, if we let $c=\max(32,2/\e,\delta^{-1})$.
\end{proof}

\providecommand{\bysame}{\leavevmode\hbox to3em{\hrulefill}\thinspace}

\end{document}